\documentclass[12pt]{amsart}

\usepackage{amsmath,amssymb,epsfig,amscd,xy,amsthm}
\xyoption{all}

\newtheorem{theorem}{Theorem}[section]
\newtheorem{proposition}[theorem]{Proposition}

\newtheorem{defn}[theorem]{Definition}

\theoremstyle{plain}

\newtheorem*{main}{Theorem 1}
\numberwithin{equation}{theorem}

\theoremstyle{remark}

\newcommand{\Kb}{{\overline K}}

\newcommand{\kbv}{{\overline k_v}}

\newcommand{\Pb}{{\bar P}}
\newcommand{\ba}{{\bar \alpha}}
\newcommand{\Qb}{{\bar Q}}

\DeclareMathOperator{\Hom}{Hom}
\DeclareMathOperator{\SL}{SL}
\DeclareMathOperator{\GL}{GL}
\DeclareMathOperator{\PGL}{PGL}
\DeclareMathOperator{\Aut}{Aut}

\DeclareMathOperator{\Spec}{Spec} 
\DeclareMathOperator{\Proj}{Proj}

\DeclareMathOperator{\Gal}{Gal}

\DeclareMathOperator{\Supp}{Supp}

\newcommand{\bP}{{\mathbb P}}

\newcommand{\lra}{\longrightarrow}
\newcommand{\ra}{\rightarrow}
\newcommand{\CO}{{\mathcal O}}
\newcommand{\fo}{{\frak o}}
\newcommand{\cI}{\mathcal{I}}

\newcommand{\fm}{{\frak m}}
\newcommand{\cU}{\mathcal{U}}
\newcommand{\cW}{\mathcal{W}}

\newcommand{\cY}{\mathcal{Y}}

\newcommand{\fp}{{\frak p}}

\newcommand{\bK}{\overline{K}}

\title{A Shafarevich-Faltings Theorem For Rational Functions}
\author{Lucien Szpiro and Thomas J. Tucker}

\begin{document}






\begin{abstract}
  Using an alternative notion of good reduction, a generalization of
  the Shafarevich finiteness theorem for elliptic curves is proved for
  self-maps of the projective line over number fields.
\end{abstract}
\maketitle

\begin{center}
  {\it This paper is dedicated to Fedor Bogomolov on the occasion of
    his sixtieth birthday.  }
\end{center}

\section{Introduction}

In 1963, Shafarevich (\cite{Sha}) proved that for any finite set $S$
of primes in a number field $K$, there are finitely many isomorphism
classes of elliptic curves with good reduction at all primes outside
of $S$ (see also \cite[IV.1.4]{Serre} for a quick proof of this
result).  Shafarevich also conjectured that there were only finitely
many isomorphism classes of abelian varieties over $K$ of any fixed
dimension, with polarization of fixed degree $d$, that have good
reduction at all primes outside of $S$.  This conjecture, often called
the Shafarevich conjecture, was later proved by Faltings
(\cite{Faltings}) as part of his celebrated proof of the Mordell
conjecture.

Shafarevich's result implies that that the minimal discriminant of an
elliptic curve can be bounded in terms of the conductor of the curve.
Indeed, this follows trivially from the fact that his result implies
that there are finitely many isomorphism classes of elliptic curves
having fixed conductor.  An explicit bound for the minimal
discriminant in terms of the conductor was later proposed by the first
author (\cite{Szpiro}).  This conjecture was proved for function
fields (\cite{Hindry, Pesenti}), but it remains open in the number
field case.  Note that the $abc$-conjecture of Masser and Oesterl{\'e}
(which is closely related to this conjecture) has been proved for
function fields by Mason \cite{Mason} (in fact, it might be more
accurate to say that Mason's result helped motivate the
$abc$-conjecture) and that Arakelov, Par{\v{s}}in, and the first
author (\cite{Arakelov, Parsin, LS2}) proved the Shafarevich
conjecture over function fields several years before Faltings proved
it for number fields.  Bogomolov, Katzarkov, and Pantev (\cite{Bogo})
have even managed to prove a version of the conjecture of
\cite{Szpiro} for hyperelliptic curves over function fields.

Perhaps the most natural definition for good reduction in the context
of nonconstant morphisms $\varphi: \bP^1_K \ra \bP^1_K$ over a number
field $K$ is to say that such a map $\varphi$ has good reduction at a
finite place $v$ if $\varphi$ extends to a map $\bP^1_{\fo_v} \ra
\bP^1_{\fo_v}$, where $\fo_v$ is the localization of the ring of
integers $\fo_K$ at $v$.  When this is the case, we will say that
$\varphi$ has {\bf simple good reduction} at $v$.  Note that this is
equivalent to saying that there is a choice of $\fo_v$-coordinates for
$\bP^1_{\fo_v}$ such that $\varphi$ can be written as $\varphi([x:y])
= [P(x,y):Q(x,y)]$ where $P$ and $Q$ are homogeneous polynomials of
the same degree in $\fo_v[x,y]$ that do not have common roots in the
residue field at $v$.  The situation here, however, is quite different
from the case of elliptic curves since any monic polynomial $f(x) \in
\fo_K[x]$ gives rise to a morphism that has simple good reduction at
all finite places.  Thus, one is led to consider alternative notions
of good reduction, one of which we will now explain.

Let $K$ be a number field, let $\fo_K$ be its ring of integers, and let
$S$ be a finite set of finite places of $\fo_K$.  We define
$$
\fo_S = \{ z \in K \quad | \quad |z|_v \leq 1 \text { for all $v
  \notin S$} \}.$$
The ring $\fo_S$ is often referred to as the ring
of $S$-integers in $K$.  Let $\bP^1_{\fo_S}$ be the usual projective
line over $\Spec \fo_S$ (which can be defined as $\Proj
\fo_S[T_0,T_1]$, as in \cite[p. 103] {H}), and let $g: \bP^1_K {\tilde
  \rightarrow} \left( \bP^1_{\fo_S} \right)_K$ be an isomorphism
(where $\left( \bP^1_{\fo_S} \right)_K$ denotes $\bP^1_{\fo_S}
\times_{\Spec \fo_S} \Spec K$, as usual).  One obtains a map $\bP^1_K
\rightarrow \bP^1_{\fo_S}$ by composing $g$ with the base extension
$\left( \bP^1_{\fo_S} \right)_K \ra \bP^1_{\fo_S}$.  For each finite
place $v \notin S$, this map gives rise to a reduction map $r_{g,v}:
\bP^1(K) \lra \bP^1(k_v)$ where $k_v$ is residue field of $\fo_K$ at
$v$.  By extending the place $v$ to a place on $\bK$, one may extend
$r_{g,v}$ to a map $r_{g,v}: \bP^1(\Kb) \lra \bP^1(\kbv)$.  This
allows us to make the following definition.

\begin{defn}
  Let $\varphi:\bP_K^1 \lra \bP^1_K$ be a nonconstant morphism of
  degree greater than 1, let $g: \bP^1_K {\tilde \rightarrow} \left(
    \bP^1_{\fo_S} \right)_K$ be an isomorphism, let $R_\varphi$ denote
  the ramification divisor of $\varphi$ over $\bK$, and let $v \notin
  S$ be a finite place of $K$ that has been extended to $\Kb$.  We say
  that $\varphi$ has {\bf critically good reduction} at a finite place
  $v \notin S$ if the following conditions are met:
\begin{enumerate}
\item for any points $P \not= Q$ in $\bP^1(\Kb)$ contained in $\Supp
  R_\varphi$, we have $r_{g,v}(P) \not= r_{g,v}(Q)$; and
\item for any points $P \not= Q$ in $\bP^1(\Kb)$ contained in
  $\varphi(\Supp R_\varphi)$ we have $r_{g,v}(P) \not= r_{g,v}(Q)$.
\end{enumerate}
\end{defn}
The terminology ``critically good reduction'' was suggested to the
authors by Joseph Silverman during the preparation of this paper.  Our
definition does not depend on how we choose to extend $v$ to all of
$\Kb$, since if $v'$ is a place of $\Kb$ that agrees with $v$ on $K$,
then there is an automorphism $\tau \in \Gal(\Kb/K)$ such that
$r_{g,v}(P) = r_{g,v'}(\tau P)$ for all $P \in \bP^1(\bK)$.  A simple
way of describing this definition is to say that all the distinct
$\bK$-points in $\Supp R_\varphi$ remain distinct after reduction at
$v$ and all the distinct $\bK$-points in $\varphi(\Supp R_\varphi)$
remain distinct after reduction at $v$.

The automorphism group $\Aut(\bP^1_{\fo_S})$ is isomorphic to
$\PGL_2(\fo_S)$ (that is, $\GL_2(\fo_S)$ modulo homotheties in
$\GL_2(\fo_S)$).  This can be seen by choosing $\fo_S$-coordinates for
$\bP^1_{\fo_S}$.  We say that two morphisms $\varphi: \bP^1_K
\rightarrow \bP^1_K$ and $\psi: \bP^1_K \rightarrow \bP_K^1$ are
$g$-equivalent if they are the same up to multiplication by an element
of $\Aut(\bP^1_{\fo_S})$ on both sides; that is, if there are
automorphisms $\gamma, \sigma \in \Aut(\bP^1_{\fo_S})$ such that
$$
\psi = \sigma_K \varphi \gamma_K,$$
where $\sigma_K$ and $\gamma_K$
are the pull-backs of $\sigma$ and $\gamma$ to $\bP^1_K$ (via the
isomorphism $g: \bP^1_K {\tilde \rightarrow} \left( \bP^1_{\fo_S}
\right)_K$ and the map $\left( \bP^1_{\fo_S} \right)_K \lra \bP^1_{\fo_S}$).  With
this terminology, the main result of this paper is the following
analog of Shafarevich's theorem (\cite{Sha}) for elliptic curves.

\begin{main}
  Let $S$ be a finite set of finite places of a number field $K$, let $n$ be
  an integer greater than one, and let $g: \bP^1_K {\tilde
    \rightarrow} \left( \bP^1_{\fo_S} \right)_K$ be an isomorphism.
  There are finitely many $g$-equivalence classes of nonconstant
  morphisms $\varphi:\bP_K^1 \lra \bP_K^1$ of degree $n$ that ramify at
  three or more points and have critically good reduction at all
  finite places outside of $S$.
\end{main}

Note that if a map $\varphi$ ramifies at exactly two points, it is
easy to see that there are automorphisms $\sigma_L$ and $\gamma_L$,
each defined over a quadratic extension $L$ of $K$, such that
$\sigma_L \varphi \gamma_L $ is a map of the form $x \mapsto a x^m$.
Thus, the equivalence classes for such maps are even easier to
describe, provided that one considers a slightly larger automorphism
group.

The proof of Theorem 1 is as follows.  In Section \ref{one}, we use a
result of Birch-Merriman and Evertse-Gy{\H{o}}ry (\cite{Birch, EG1})
to show that there is a finite set $\cY$ such that $\Supp R_\varphi$
and $\varphi(\Supp R_\varphi)$ are both contained in $\cY$ after the
application of suitable automorphisms.  Then, in Section \ref{two}, we
apply a result of Mori (\cite{Mori}) to conclude that this gives us a
finite set of equivalence classes of maps.\\
\\
\noindent{\it Acknowledgments.}  We would like to thank Joseph Silverman for
many helpful conversations.  We would like to thank Xander Faber for
his careful reading of the paper and for his many useful
suggestions.

\section{Connections with other notions of good reduction}
Note that the notion of critically good reduction is in some ways
quite similar to the notion of good reduction on an elliptic curve.  A
model $y^2 = F(x)$ for an elliptic curve $E$ over $\fo_K$ has good
reduction at a finite place $v \nmid 2$ if and only if $F$ has
distinct roots and the leading coefficient of $F$ is a $v$-adic unit.
This implies that the ramification points of the map obtained by
projecting onto the $x$-coordinate remain distinct after reduction at
$v$.  In fact, the multiplication-by-two map on an elliptic curve is
associated to a map $\varphi: \bP_K^1 \lra \bP_K^1$, called a
Latt{\`e}s map, which can be written explicitly as $\varphi(x) =
\frac{(F'(x))^2 - 8 x F(x)}{4 F(x)}$.  This map is simply the usual
rational function that describes the $x$-coordinate of $2\beta$ in
terms of the $x$-coordinate of $\beta$ for $\beta$ a point on $E$ (see
\cite[Chapter 2]{AEC}, for example).  Note that our definition of
critically good reduction depends on our choice of an isomorphism
$g: \bP^1_K {\tilde \rightarrow} \left( \bP^1_{\fo_S} \right)_K$.
Choosing coordinates $[x:y]$ gives rise to a natural isomorphism $g$,
and when a map $\varphi$ is written explicitly in terms of coordinates
we will say that $\varphi$ has critically good reduction at $v$ if it
has critically good reduction at $v$ in the model these coordinates
determine.  With this convention, we have the following.

\begin{proposition}\label{good good}
  Let $E$ be an elliptic curve over a number field $K$ and let $S$ be
  a set of finite places of $K$ containing all places $v$ such that $v
  | 2$.  If the model $y^2 = F(x)$ for $E$ over $\fo_S$ has good
  reduction at all finite places $v \notin S$, then the corresponding
  Latt{\`e}s map $\varphi(x) = \frac{(F'(x))^2 - 8 x F(x)}{4 F(x)}$
  has both critically good reduction and simple good reduction at all
  finite places $v \notin S$.
\end{proposition}
\begin{proof}
  Let $v \notin S$ be a finite place of $K$ and let $\fm$ denote the
  maximal ideal corresponding to $v$ in $\fo_S$.  Since the model $y^2
  = F(x)$ has good reduction at $v$, the leading coefficient of $F$ is
  a unit at $v$ and the roots of $F$ are distinct at $v$.  Suppose
  that some $\alpha$ in the algebraic closure of $\fo_S / \fm$ is a
  root of both $4F(x)$ and $(F'(x))^2 - 8 x F(x)$ modulo $\fm$; then
  it is a multiple root of $F(x)$ modulo $\fm$, which would mean that
  $F$does not have distinct roots modulo $\fm$ and that $E$ therefore
  not have good reduction at $v$.  Thus, there is no such $\alpha$, so
  $\varphi(x)$ is well-defined modulo $\fm$ for all $x$.  Since the
  leading coefficient of $(F'(x))^2 - 8 x F(x)$ is the same as the
  leading coefficient of $F$ (and is thus a unit at $v$), the map
  $\varphi$ is also well-defined at infinity.  Hence, $\varphi$ has
  simple good reduction at $v$.
  
  To see that the map $\varphi(x) = \frac{(F'(x))^2 - 8 x F(x)}{4
    F(x)}$ has critically good reduction, we simply note that the
  ramification points of this map are the $x$-coordinates of the the
  points in $E[4] \setminus E[2]$ (i.e., the $4$-torsion points of $E$
  that are not $2$-torsion points) and that their images are the
  $x$-coordinates of the points $E[2] \setminus \{ 0 \}$. The
  reduction map at $v$ is injective on prime-to-$p$ torsion for $v
  \mid p$ (see \cite[Proposition VII.3.1]{AEC}, for example), and $p
  \not= 2$ by assumption, so all of the points in $E[4]$ are distinct
  modulo $\fm$.  Thus, all of the points in $\Supp R_\varphi$ and
  $\varphi(\Supp R_\varphi)$ are distinct modulo $\fm$, as desired.
\end{proof}

More generally, it is possible to have simple good reduction without
having critically good reduction.  This can be seen, for example, by
taking any monic polynomial $f(x) \in \fo_K[x]$ such that $f'(x)$ has
multiple roots at some place $v$.  It is also possible to have
critically good reduction without having simple good reduction; the
map $x \mapsto p x^2$ is an example of this.  On the other hand, under
fairly generic hypotheses, critically good reduction does imply simple
good reduction.  

\begin{proposition}
  Let $\varphi(x) = P(x)/Q(x)$ be a rational function of degree $d$
  with coefficients in $\fo_S$ for $S$ some finite set of finite
  places of a number field $K$.  Let $v \notin S$ be a finite place of
  $K$.  Suppose that $\varphi$ has $2d - 2$ distinct ramification
  points and that the leading coefficients of $P$, $Q$, and $P'(x)Q(x)
  - P(x) Q'(x)$ are all $v$-adic units.  Then, if $\varphi$ has
  critically good reduction at $v$, it also has simple good reduction
  at $v$.
\end{proposition}
\begin{proof}
  Suppose that $\varphi$ does not have simple good reduction.  Then
  there is a $\alpha$ with $|\alpha|_v \leq 1$ such that $P(\alpha)$
  and $Q(\alpha)$ are both zero modulo the maximal ideal at $v$.  Let
  $\Pb$, $\Qb$, and $\ba$ denote the reductions of $P$, $Q$, and
  $\alpha$, respectively, at $v$.  Then $\Pb'(\ba) \Qb(\ba) - \Pb(\ba)
  \Qb'(\ba) = 0$; taking the derivative again, we find that
  $\Pb''(\ba) \Qb(\ba) - \Pb(\ba) \Qb''(\ba) = 0$.  Thus, $\ba$ is a
  double root of $\Pb' \Qb - \Pb \Qb'$.  It follows that $\Pb' \Qb -
  \Pb \Qb'$ has fewer distinct roots than $P'Q - PQ'$.  Since $'Q -
  PQ'$ and $\Pb' \Qb - \Pb \Qb'$ both have the same degree, it follows
  that two roots of $P'Q - PQ'$ reduce to the same root of $\Pb' \Qb -
  \Pb \Qb'$.  Since the roots of $P'Q - PQ'$ are all ramification
  points of $\varphi$, this means that $\varphi$ does not have
  critically good reduction at $v$.
\end{proof}

Shafarevich's theorem (\cite{Sha}) can be considered a special case of
Theorem 1 by taking the map $\psi_E: \bP_K^1 \lra \bP_K^1$
corresponding to the multiplication-by-4 map on an elliptic curve $E$
(this is simply $\varphi \circ \varphi$ where $\varphi$ is the
Latt{\`e}s map considered earlier).  As noted in \cite[IV.1.4]{Serre},
if $E$ has good reduction outside of a finite set of places $S$ that
includes the places lying over 2 and for which $\fo_S$ is a principal
ideal domain, then there is a model for $E$ given by an equation $y^2
= F(x)$, where the coefficients of $F$ are in $\fo_S$, the leading
coefficient of $F$ is an $S$-unit, and $F$ does not have multiple
roots at any of the finite places outside of $S$.  The map $\psi_E$
ramifies over the points where the projection onto the $x$-coordinate
map from $E$ to $\bP^1$ ramifies (note that the Latt{\`e}s map on
$\bP^1$ corresponding to multiplication-by-2 fails to ramify at
infinity, which is why use multiplication-by-4) and has critically
good reduction outside of $S$, since it ramifies at the
$x$-coordinates of the points in $E[8] \setminus E[2]$ and over the
$x$-coordinates of the points in $E[2]$.  Now, let $E'$ be another
elliptic curve with good reduction outside $S$ having the property
that $\psi_{E'}$ is $g$-equivalent to $\psi_E$.  Let $y^2 = G(x)$ be a
model for $E'$, where $G$ has coefficients in $\fo_S$ and the leading
coefficient of $G$ is an $S$-unit.  Since $\psi_E$ and $\psi_{E'}$ are
$g$-equivalent, we have an automorphism $\gamma \in
\Aut(\bP^1_{\fo_S})$ that takes the $x$-coordinates of the points in
$E[2]$ to the $x$-coordinates of the points in $E'[2]$.  After this
automorphism we may suppose that $F$ and $G$ have the same roots
(since their roots, along with $\infty$ are the $x$-coordinates of the
2-torsion).  This means that $G = u F$ for some $S$-unit $u$.  Since
replacing $y$ with a multiple of $y$ in the Weierstrass equation for
an elliptic curve does not change the isomorphism class, we see that
the isomorphism class of $E'$ is determined is determined by the coset
class of $u$ in $\fo_S^* / (\fo_S^*)^2$.  Since $\fo_S^*$ is finitely
generated, it follows that there are only finitely many isomorphism
classes of elliptic curves $E'$ such that $\psi_{E'}$ is
$g$-equivalent to $\psi_E$.  Hence, Theorem 1 implies the Shafarevich's
theorem.  Note that Shafarevich's theorem can also be obtained
directly from the results of Birch-Merriman and Evertse-Gy{\H{o}}ry
(\cite{Birch, EG1}).

\section{Finiteness theorems and homogeneous forms}\label{one}
Let $K$ be number field, let $S$ be a finite set of finite places of $K$,
and let $g: \bP^1 \lra \left(\bP_{\fo_S} \right)_K$ be an isomorphism.
Let $v \notin S$ be a finite place of $K$.  Let $k_v$ denote the
residue field of $K$ at $v$.  The isomorphism $g$ gives us a reduction
map
$$
r_{g,v}: \bP^1(K) \lra \bP^1(k_v),$$
as described in the introduction.  We
may extend $v$ to all of $\Kb$.  We then have a reduction map
$$ r_{g,v}: \bP^1(\Kb) \lra \bP^1(\kbv).$$ 

\begin{defn}
  Let 
  $$\cU = \{ z_1, \dots, z_n \},$$
  where $z_i$ are distinct elements
  of $\bP^1(\Kb)$, and let $v \notin S$ be a finite place of $K$.  We
  say that the set $\cU$ is $r_{g,v}$-distinct if
  $r_{g,v}(z_i) \not= r_{g,v}(z_j)$ for all $i \not= j$.
\end{defn}

Using coordinates, we have a very simple criterion for deciding when a
set is $r_{g,v}$-distinct.  We may choose coordinates for $\bP^1_K$
and $\bP^1_{\fo_S}$ such that the map $\bP^1_K \ra \bP^1_{\fo_S}$
obtained by composing $g$ with the base extension map is given in
coordinates by $[a:b] \mapsto [ca :cb]$ for any $c \in K$ such that
$ca, cb \in \fo_S$.  For each $z_i \in \cU$, we write $z_i$ as
$[a_i:b_i]$ where $a_i, b_i \in \Kb$ and $\max(|a_i|_v, |b_i|_v) = 1$.
Then $\cU$ is $r_{g,v}$-distinct if and only if
\begin{equation}\label{real easy}
|a_i b_j - a_j b_i|_v = 1
\end{equation}
for all $i \not= j$.  This follows from the fact that 
$$
r_{g,v}([a_i: b_i]) = [a_i \mod \fm_v: b_i \mod \fm_v] \in
\bP^1(\kbv)$$
where $\fm_v$ is the maximal ideal corresponding to $v$
in the ring of integers of $\bK$.

\begin{defn}
  Let $S$ be a finite set of finite places of $K$.  We say that a set $\cU
  \subseteq \bP^1 (\Kb)$ is $S$-good if for every finite place $v \notin S$,
  the set $\cU$ is $r_{g,v}$-distinct.
\end{defn}

We will relate the notion of sets being $S$-good to the discriminant
of homogeneous forms that vanish at $\cU$.  This will allow us to
apply a finiteness result due to Birch and Merriman (\cite{Birch};
see also \cite{EG1}).  We define

$$
\fo^*_S = \{ z \in K \quad | \quad |z|_v = 1 \text { for all $v
  \notin S$} \}.$$
The elements of this unit group $\fo_S^*$ are
called {\it $S$-units}.  For any $\delta \in \fo_K$ we define $\delta
\fo_S^*$ to be the set of all elements of $K$ of the form $\delta u$
for $u \in \fo_S^*$.

Let $F(x,y)$ be a homogeneous form of degree $n$ in $\fo_S[x,y]$.
Factoring $F$ in $\Kb[x,y]$ as
$$
F(x,y) = \prod_{i=1}^n (\beta_i x - \alpha_i y),$$
we define
the discriminant $\Delta(F)$ of $F$ as
$$ \Delta(F) = \left( \prod_{i < j} (\alpha_i \beta_j - \beta_i
  \alpha_j) \right)^2.$$

For $\gamma \in
\SL_2(\fo_S)$, with 
$$
\gamma = 
\left( \begin{array}{rr}
    a & b \\
    c & d \end{array} \right),$$
we let
$$
\gamma(F(x,y)) = F(ax + by, cx + dy).$$
Let $F(x,y)$ and $G(x,y)$
be homogeneous forms in $\fo_S[x,y]$.  We say that $F$ and $G$ are in
the same $\fo_S$-equivalence class if there is an element $\gamma \in
\SL_2(\fo_S)$ and a $\lambda \in \fo_S^*$ such that $\gamma(F(x,y)) =
\lambda G(x,y)$.  Birch-Merriman (\cite[Theorem 1]{Birch}) and
Evertse-Gy{\H{o}}ry (\cite[Theorem 3]{EG1}) proved the following
theorem about $\fo_S$-equivalence classes of homogeneous forms.

\begin{theorem}\label{BM} (\cite[Theorem 1]{Birch}, \cite[Theorem 3]{EG1}.)  Let
  $\delta$ be any nonzero element of $\fo_S$ and let $n$ be a positive
  integer.  There are finitely many $\fo_S$-equivalence classes of
  degree $n$ homogeneous forms $F(x,y) \in \fo_S[x,y]$ such that
  $\Delta(F) \in \delta \fo_S^*$.
\end{theorem}

We note that Birch and Merriman only state their result for forms in
$\fo_K[x,y]$ and in the case that $\delta = 1$, but one can deduce the
statement above for forms in $\fo_S[x,y]$ and arbitrary $\delta \in
\fo_K$ from their result (\cite[Theorem 1]{Birch}).  Evertse and
Gy{\H{o}}ry's result (\cite[Theorem 3]{EG1}) is an effective version
of Theorem \ref{BM}.
 
\begin{proposition}\label{from BM}
  Let $n$ be an integer.  Then there is a finite set $\cY$ such that
  for any $S$-good $\Gal({\overline K}/K)$-stable set $\cU$ of
  cardinality $n$, there is a $\gamma \in \Aut(\bP^1_{\fo_S})$ such
  that $\gamma_K(\cU) \subseteq \cY$.
\end{proposition}
\begin{proof}
  We will need a little notation to deal with the fact that $\fo_S$
  may not be a unique factorization domain.  For any fractional ideal
  $J$ of $\fo_S$ we let $v(J) = e_{J_v}$ where $e_{J_v}$ is the power
  of the prime $\fp_v$ corresponding to $v$ in the factorization of
  $J$ into prime ideals.  Let $I_1, \dots, I_s$ be a set of (integral)
  ideals in $\fo_S$ representing the ideal classes of $\fo_S$; that
  is, for any fractional ideal $J$ of $\fo_S$, there is an $\alpha \in
  K$ such that $\alpha J = I_\ell$ for some $\ell$.  Then there are
  finitely many sets of the form $\delta \fo^*_S$ where $\delta$ is an
  element of $\fo_k$ such that $v((\delta)) \leq (2n - 2)v(I_j)$ for all
  $I_j$.  Let $\cW = \{ \delta_1 \fo_S^*, \dots, \delta_m \fo_S^* \}$
  be the set of all such $\delta \fo_S^*$.
  
  By Theorem \ref{BM}, for any $\delta \in \fo_S$, there are at most
  finitely many $\fo_S$-equivalence classes of forms in $G \in
  \fo_S[x,y]$ of degree $n$ such that $\Delta(G) \in \beta \fo_S^*$.
  Thus we may choose a set of forms $G_1, \dots, G_t$ such that for
  any form $F$ of degree $n$ with $\Delta(F) \in \delta_j \fo_S^*$,
  for $\delta _j \fo^* \in \cW$, there is some $G_i$ such that $F$ is
  $\fo_S$-equivalent to $G_i$.  Let
  $$ \cY = \{ [a:b] \in \bP^1(\Kb) \; \mid \; G_i(a,b) = 0 \text{ for some
    $G_i$} \}.$$

  Let $\cU = \{ [a_1:b_1], \dots, [a_n:b_n] \}$ be $S$-good and let
  $H$ be any homogeneous form of degree $n$ in $\fo_S[x,y]$ that
  vanishes on $\cU$.  After multiplying through by a nonzero element
  $\alpha \in K$, we obtain a form $\alpha H \in \fo_S[x,y]$ such that
  the coefficients of $\alpha H$ generate one of the ideals $I_j$
  above.  Write each $[a_i:b_i]$ as $[a'_i:b'_i]$ where
  $\max(|a'_i|_v, |b'_i|_v) = 1$.  Then, there is an element $\kappa
  \in {\overline K}$ such that
  $$
  \alpha H(x,y) = \kappa \prod_{i=1}^n(b'_i x - a'_i y).$$
  We have
  $$
  | \Delta(\alpha H)|_v = |\kappa|_v^{2n-2} \left| \left( \prod_{i
        < j} (a'_i b'_j - b'_i a'_j) \right)^2 \right|_v =
  |\kappa|^{2n -2}_v $$
  since $|a'_i b'_j - b'_i a'_j|_v = 1$ for all
  $i,j$ (because $\cU$ is $S$-good) and multiplying a form of degree
  $d$ through by a constant $\kappa$ changes the discriminant by a
  factor of $\kappa^{2n - 2}$.  Now, $v((\kappa)) \geq 0$ by the Gauss
  lemma for polynomials and $v((\kappa)) \leq v(I_j)$ because of the
  fact that the coefficients of $\alpha H$ generate $I_j$.  Since
  $v((\kappa)) \leq v(I_j)$, we see that $\kappa^{2n - 2} \fo_S^* = \delta_i
  \fo_S^*$ for some $\delta_i \fo_S^* \in \cW$.
   
  Thus, for some $G_\ell$, there is a $\tau \in \SL_2(\fo_S)$ and
  $\lambda \in \fo_S^*$ such that $\tau(\alpha H) = \lambda G_\ell$.  Hence,
  for each $[a_i:b_i]$, we have $G_\ell(\tau([a_i:b_i])) = 0$.  Our
  choice of coordinates gives an inclusion of $\SL_2(\fo_S)$ into
  $\Aut(\bP^1_{\fo_S})$.  Thus, $\tau$ corresponds to an element
  $\gamma \in \Aut(\bP^1_{\fo_S})$ such that
  $$
  \gamma_K(\cU) \subseteq \cY,$$
  as desired.
\end{proof}

\section{Applying a Result of Mori}\label{two}

We will now state a result due to Mori (\cite{Mori}).  We note that
what Mori proves is much more general than what is required here.  Let
$A$ and $B$ be schemes of finite type over an algebraically closed
field $L$, and let $Z$ be a closed subscheme of $A$.  Let $p: Z \lra
B$ be an $L$-morphism.  Let $\Hom_L(A,B;p)$ be the set of $L$-morphisms
from $A$ to $B$ that extend $p$, that is
$$
\Hom_L(A,B;p) = \{ f: A \lra B \; \mid \; \text{$f$ is an
  $L$-morphism and $f |_Z = p$} \}.$$
We let $\cI_Z$ denote the ideal
sheaf of $Z$ in $A$ and let $T_B$ denote the tangent sheaf of $B$ over
$L$.

\begin{theorem}\label{Mori thm} (\cite[Propositions 1 and 3]{Mori}.)
  The set $\Hom_L(A,B;p)$ is represented by a closed subscheme of
  $\Hom_L(A,B)$ and for any closed point $f$ in $\Hom_L(A,B;p)$, we
  have
  $$
  T_{f, \Hom_L(A,B;p)} \cong H^0(A, f^*T_B \otimes_{\CO_A}
  \cI_Z),$$
  where $ T_{f, \Hom_L(A,B;p)}$ is the tangent space of $f$
  in $\Hom_L(A,B;p)$ over $L$.
\end{theorem}

Using this theorem, we are able to derive the following proposition.  

\begin{proposition} \label{from Mori}
  Let $\cY$ be a finite subset of $\bP^1(L)$ and let $n > 1$ be an
  integer.  Then there are finitely many morphisms $\varphi: \bP^1_L
  \lra \bP^1_L$ of degree $n$ satisfying all of the following
  conditions:
\begin{enumerate}
\item $\Supp
  R_\varphi \subseteq \cY$;
\item $\varphi(\Supp R_{\varphi}) \subseteq \cY$; and
\item $|\Supp R_{\varphi}| \geq 3$.
\end{enumerate}
\end{proposition}
\begin{proof}
  There are at most $|\cY|^n$ possible divisors $D$ with support in
  $\cY$ that could be ramification divisors of morphisms
  $\varphi:\bP^1_L \rightarrow \bP^1_L$ of degree $n$.  Furthermore,
  the fact that $\cY$ is finite means that there are finitely many
  possibilities for the image $\varphi(\Supp R_\varphi)$.  Thus, it
  suffices to show that for any divisor $D = \sum_{i=1}^m f_i Q_i$ on
  $\bP^1$ with $|\Supp D| \geq 3$ and any sequence of points $P_1,
  \dots, P_m$ (not necessarily distinct) with $P_i \in \cY$, there are
  at most finitely many morphisms $\varphi$ of degree $n$ such that
  $R_{\varphi} = D$ and $\varphi(Q_i) = P_i$.

  Let $D$ be the ramification divisor of a map $\varphi$ of degree
  $n$.  We write $D = \sum_{i=1}^m (e_{Q_i} - 1) Q_i$ where $e_{Q_i}$
  is the ramification index of $\varphi$ at ${Q_i}$. For each ${Q_i}
  \in \Supp D$, let $\cI_{Q_i}$ represent the ideal sheaf of ${Q_i}$
  in $\bP^1$, and let $Z$ be the subscheme of $\bP^1$ with ideal sheaf
  $\cI = \prod_{i=1}^m \cI^{e_{Q_i}}$ (note that the correct exponent
  here is $e_{Q_i}$, not $(e_{Q_i} - 1)$).

  At each $P_i$ we have
  $$\varphi^*(\cI_{P_i}) = \prod_{\varphi(Q) = P_i} \cI_Q^{e_Q}$$
  (by
  the definition of the ramification index.) Thus, at each $Q$, the
  map to $P_i \in \bP_L^1$ from the scheme defined by
  $\cI_{Q_i}^{e_{Q_i}}$ is induced by the unique nonzero map of
  $L$-algebras from $L$ to $L[x]/x^{e_Q}$.  These piece together to
  form a unique map $p: Z \lra \bP_L^1$.  Thus, we have $R_\varphi =
  D$ and $\varphi(Q_i) = P_i$ for $i=1, \dots, m$ exactly when
  $\varphi$ restricts to $p$ on $Z$.

   Using the Riemann-Hurwitz formula, we see that
   \begin{equation*}
   - \deg \cI_Z = \deg R_\varphi + |\Supp D| = (2n - 2) + |\Supp D| \geq
   2n +1.
   \end{equation*}
   Since $\deg \varphi^* T_{\bP_L^1} = 2n$, we have $\deg (\varphi^* T_{\bP_L^1}
   \otimes \cI_Z) < 0$, so
   $$
   \dim_L T_{\varphi, \Hom_L(\bP_L^1, \bP_L^1;p)} = \dim_L H^0(\varphi^*
   T_{\bP_L^1} \otimes \cI_Z) = 0. $$
   The scheme
   $\Hom_L(\bP_L^1, \bP_L^1;p)$ therefore has dimension zero.  Since
   it is also Noetherian, this means that it is finite.  Thus, there
   are at most finitely many maps $\varphi$ such that $R_\varphi = D$
   and $\varphi(Q_i) = P_i$ for $i=1, \dots, m$.  This completes our
   proof.
\end{proof}

We are now ready to prove Theorem 1.

\begin{proof} ({\it Of Theorem 1.})
  Since $\Supp R_\varphi$ and $\Supp (\varphi_* R_\varphi)$ are both
  $S$-good, there are $\gamma, \sigma \in \SL_2(\fo_S) \subseteq
  \Aut(\bP^1_{\fo_S})$ such that both $\gamma_K (\Supp R_\varphi)$ and
  $\sigma_K (\varphi(\Supp R_\varphi))$ are contained in $\cY$.  Then
  the map
$$ \psi = \sigma_K \varphi \gamma_K^{-1}$$
satisfies the conditions of Proposition \ref{from Mori}.
\end{proof}

\section{Generalizations of Theorem 1}
There are many possible ways in which Theorem 1 might be generalized
and strengthened.  

\subsection{Effectivity}
As noted earlier, Evertse and Gy{\H{o}}ry (\cite{EG1}) have proved an
effective version of Theorem \ref{BM}.  More precisely, they are able
to produce an explicit constant $C$ (depending only on $S$ and the
degree $n$) such that each $\fo_S$-equivalence class of homogeneous
forms contains a form with height less than $C$.  This translates
immediately into a bound on the height of points in the set $\cY$ used
in Proposition \ref{from BM}.  While Theorem \ref{Mori thm} is not
effective as stated, it should be possible to derive an effective
version of Proposition \ref{from Mori} by viewing the conditions
placed on $\varphi$ as hypersurfaces in a suitable space of rational
functions and applying an arithmetic B{\'e}zout-type theorem such as
the one proved by Bost, Gillet, and Soul{\'e} in \cite{BGS}.  We plan
to treat this question in a later paper.  Note that effective versions
of the Shafarevich conjecture have been proved in the number field
case for elliptic curves by Silverman and Brumer (\cite{SilverBru})
and for more general curves in the function field case by Caporaso and
Heier (\cite{Caporaso, GH}).

\subsection{Higher dimensions}
It should be possible to formulate a higher-dimensional version of
Theorem 1.  One might for example impose the condition that all the
components of the ramification divisor and its image intersect
properly at all primes outside of a finite set $S$.  Mori's results
still apply in higher dimensions.  What is less clear is how to apply
finiteness results about homogeneous forms.

\subsection{Function fields}

It may also be possible to prove an analog of Theorem 1 for maps over
function fields over finite fields.  Note, however, that since there
are maps in characteristic $p$ that ramify over single point, certain
classes of maps would probably have to be excluded.  Unfortunately,
these classes include the maps that corresond to Drinfeld modules.
Taguchi (\cite{Taguchi}) has also shown that the Shafarevich
conjecture with the usual notion of good reduction does not hold for
Drinfeld modules.  It would be interesting to find a notion of good
reduction that gives rise to an analog of the Shafarevich conjecture
that {\it does} hold for Drinfeld modules.

\subsection{An $abc$-conjecture for morphisms of the projective line}
In analogy with the $abc$-conjecture and the conjecture of
\cite{Szpiro}, one might define a ``critical conductor'' for a
morphism $\varphi:\bP^1_K \ra \bP^1_K$ and ask if there is some analog
of the minimal discriminant of an elliptic curve that might be bounded
explicitly in terms of this critical conductor.  One candidate for
this analog of the minimal discriminant would be some sort of minimal
resultant for polynomials $P$ and $Q$ where $\varphi$ is written as
$P/Q$.  Since the resultant of $P$ and $Q$ will vary over
$g$-equivalence classes of morphisms, it is not clear what sort of
dependence one might expect.


\def\cprime{$'$}

\end{document}